\documentstyle[amsfonts,amssymb,amsmath,12pt]{article}
\newtheorem{thm}{Theorem}

\newtheorem{lem}[thm]{Lemma}

\date{}
\begin{document}
\title
{ Representations of the braid group $B_n$ and the highest weight
modules of $U(\mathfrak{sl}_{n-1})$ and $U_q(\mathfrak{sl}_{n-1})$
}
\author{Alexandre V. Kosyak, Inst.
of Math.
Kiev/MPI, Bonn\\
kosyak01@yahoo.com, kosyak@imath.kiev.ua
\footnote{The author would like to thank the Max-Planck-Institute
of Mathematics and the Institute of Applied Mathematics,
University of Bonn for the hospitality. The partial financial
support by the DFG project 436 UKR 113/87 is gratefully
acknowledged.}} \maketitle

\begin{abstract}
In \cite{KosAlb07q} we have constructed a
$\left[\frac{n+1}{2}\right]+1$ parameters family of irreducible
representations of the Braid group $B_3$ in arbitrary dimension
$n\in {\mathbb N}$, using a $q-$deformation of the Pascal
triangle. This construction extends in particular results by
S.P.~Humphries (2000), who constructed  representations of the
braid group $B_3$ in arbitrary dimension using the classical
Pascal triangle. E.~Ferrand (2000) obtained  an equivalent
representation of $B_3$ by considering two special operators in
the space ${\mathbb C}^n[X].$ Slightly more general
representations were given by I.~Tuba and H.~Wenzl (2001). They
involve $[\frac{n+1}{2}]$ parameters (and also use the classical
Pascal's triangle). The latter authors also gave  the complete
classification of all simple representations of $B_3$ for
dimension $n\leq 5$.  Our construction generalize all mentioned
results and throws a new light on  some of them. We also study the
irreducibility and equivalence of the constructed representations.

 In the present article we show that all representations constructed
in \cite{KosAlb07q} may be obtained by taking  exponent  of the
highest weight modules of $U(\mathfrak{sl}_2)$ and
$U_q(\mathfrak{sl}_2)$. {\it We  generalize these connections}
between the representation of  the braid group $B_n$ and  the
highest weight modules of the $U_q({\mathfrak sl}_{n-1})$ {\it for
arbitrary} $n$ using the well-known {\it reduced Burau
representations}.
\end{abstract}

\newpage
\section{Introduction. Braid group representations }
Our {\it aim} is to describe the {\it dual} $\hat{B_n}$ of the
{\it braid group} $B_n$. It is natural to compare the {\it
representation theory} of the {\it symmetric group}  $S_n$ and of
the braid group $B_n$. We know almost everything about
representation theory of the symmetric group $S_n$. We know the
description of the {\it dual} $\hat{S_n}$ in terms of  {\it Young
diagrams}. We know even the {\it Plancherel measure} on
$\hat{S_n}$. The {\it Young graph} explains how to decompose the
restriction $\pi\mid_{S_{n-1}}$ of the representation $\pi\in
\hat{S_n}$, etc.

The  braid groups $B_n$ are {\it defined} by the generators
$\sigma_i,\,1\leq i\leq n-1$ and by the relations
$\sigma_i\sigma_{i+1}\sigma_i=\sigma_{i+1}\sigma_i\sigma_{i+1},\quad
\sigma_i\sigma_j=\sigma_i\sigma_j$ for $\mid i-j\mid\,\geq 2$. The
{\it dual} $\hat{B_n}$ of the group $B_n$ {\it is known} only for
the {\it commutative case} when $n=2$. In this case
$B_2\cong{\mathbb Z}$ hence $\hat{B_2}\cong S^{1}$. The {\it
representation theory} for the braid groups $B_n$ is much more
{\it complicated} than for $S_n$. The {\it reason} is the
following. In the case of the group $S_n$ we have the essential
({\it quadratic) relation} $\sigma_i^2=1$, hence
$Sp\,(\pi(\sigma_i))\subseteq\{-1,1\}$. In the case of the group
$B_n$ we do not have these conditions. Since
$\sigma_i\sigma_{i+1}\sigma_i=\sigma_{i+1}\sigma_i\sigma_{i+1}$ we
have $Sp\,(\pi(\sigma_i))=Sp\,(\pi(\sigma_{i+1}))$, but the {\it
spectra} $Sp\,(\pi(\sigma_i))$  may be almost {\it arbitrary}.

The {\it Hecke algebra} $H_n(q)$ see f.e.\cite{Jon87} appears as
the factor algebra of the group algebra of the group $B_n$ subject
to the following  {\it quadratic relation}
$\sigma_i^2=(q-1)\sigma_i+q,\,1\leq i\leq n-1$, hence
$Sp\,(\pi(\sigma_i))\subseteq\{-1,q\}$ and $H_n(q)\cong {\mathbb
C}[S_n]$. This is a reason why the representation theory of  Hecke
algebras is well developed.

The {\it next step} is to impose the {\it polynomial condition}
$p_k(\sigma_i)=0$ on the generators $\sigma_i$ where $k$ is the
order of the polynomial $p_k(x)$. For $k=3$ the corresponding
algebra is called {\it Birman--Murakami--Wenzl type algebra} or
simple BMW algebra see \cite{Mur87,Wen98} (see also
\cite{OgiPia05} ) and so on.

The situation becomes much more complicated if no additional
conditions on the spectra are imposed. We {\it shall study} this
{\it general case} for .

In  \cite{TubWen01} I.Tuba and H.Wenzl  gave  the {\it complete
classification} of all {\it simple representations} of $B_3$ for
{\it dimension} $\leq 5$.

In \cite{ForWSV03} E.Formanek et al. gave the {\it complete
classification} of all {\it simple representations of $B_n$} for
{\it dimension} $\leq n$.

 We {\it generalize the results} I.Tuba and H.Wenzl for $B_3$,
give {\it new representations} of  $B_n$ for {\it large dimension}
and establish {\it connection} between the {\it representations}
of $B_n$ and {\it the highest weight modules} of the {\it quantum
group} $U_q({\mathfrak sl}_{n-1})$.

More precisely, in the work \cite{KosAlb07q} with S.Albeverio we
have constructed a $\left[\frac{n+1}{2}\right]+1$ parameter family
of irreducible representations of the braid group $B_3$ {it in
arbitrary dimension} $n\in {\mathbb N}$, using {\it a
$q-$deformation of the Pascal triangle}. This construction extends
in particular results by S.P.~Humphries \cite{Hum00}, I.~Tuba and
H.~Wenzl \cite{TubWen01}, and  E.~Ferrand \cite{Fer05}. The {\it
irreducibility} and the {\it equivalence} of the  constructed
representations is studied. For example the representations
corresponding to different $q$ and $n$ are {\it nonequivalent}.

In this article we show that there is a striking {\it connection}
between these {\it representations} of $B_3$ and a highest weight
modules of the {\it quantum group} $U_q({\mathfrak sl}_2)$, a
one-parameter {\it deformation of the universal enveloping
algebra} $U({\mathfrak sl}_2)$ of the Lie algebra ${\mathfrak
sl}_2$. The starting point for all these considerations is some
homomorphism $\rho_3$ of the braid group $B_3$ into ${\rm
SL}(2,{\mathbb Z}):$
$$\rho_3:B_3\mapsto {\mathfrak sl}_2 \stackrel{\exp}{\mapsto} {\rm
SL}(2,{\mathbb Z})$$
$$
\sigma_1\mapsto \left(\begin{smallmatrix}
0&1\\
0&0
\end{smallmatrix}\right)
\stackrel{\exp}{\mapsto} \left(\begin{smallmatrix}
1&1\\
0&1
\end{smallmatrix}
\right),\,\, \sigma_2\mapsto \left(\begin{smallmatrix}
0&0\\
-1&0
\end{smallmatrix}\right)
\stackrel{\exp}{\mapsto} \left(\begin{smallmatrix}
1&0\\
-1&1
\end{smallmatrix}\right).
$$

The constructed representations may be treated as the {\it
$q-$symmetric power} of this {\it fundamental representation} or
as an appropriate {\it $q-$exponen\-tial} of the
 highest weight modules of $U_q({\mathfrak sl}_2)$.

{\it We  generalize these connections} between the representation
of  the braid group $B_n$ and  the highest weight modules of the
$U_q({\mathfrak sl}_{n-1})$ {\it for arbitrary} $n$ using the
well-known {\it reduced Burau representation} $b_{n}^{(t)}$ see
c.f. \cite{Jon87}. We note that in particular
$\rho_3=b_{3}^{(-1)}$.

Let ${\mathfrak g}$ be the Lie algebra defined by a Cartan matrix
${\bf A}$ and let ${\bf B}$ be the corresponding braid group.
Denote by ${\bf U}({\mathfrak g})$ the quantized enveloping
algebra of ${\mathfrak g}$ over the field ${\mathbb C}(v)$, and
let $V$ be the integrable ${\bf U}({\mathfrak g})-$module. In
\cite{Lusz93} G.~Lusztig defined a natural action of ${\bf B}$ on
$V$ which permutes the weight space of $V$  according to the
action of the Weyl group on the weights. This rather {\it general
but different approach} allows us also to construct the
irreducible representations of the braid group ${\bf B}$ (see
\cite{OhKaKw96}).
%
\vskip 0.5cm
\newpage
{\bf 0. Definition of the Artin braid group $B_n$ }
$$
B_n=\langle(\sigma_i)_{i=1}^{n-1},
\mid\sigma_i\sigma_{i+1}\sigma_i=\sigma_{i+1}\sigma_i\sigma_{i+1},
\quad \sigma_i\sigma_j=\sigma_i\sigma_j,\quad\mid i-j\mid\,\geq
2\rangle.
$$
$B_n=\pi_1(X)$ is the {\it fundamental group} $\pi_1$ of the {\it
configuration space} $X=\{{\mathbb C}^n\setminus\Delta\}/S_n$
where $\Delta=\{(z_1,...,z_{n})\mid x_i=z_j\text{\,\,for\,\,
some\,\,}i\not=j\}$ and the group $S_n$ act freely on ${\mathbb
C}^n\setminus\Delta$ by permuting coordinates.

A {\bf BRAID} on $n$ strings {\it is a collection of curves in
${\mathbb R}^3$ joining $n$ points in a horizontal plane to the
$n$ points directly below them on another horizontal plane}.
Operation: {\it concatenation}.

{\unitlength=0.6mm
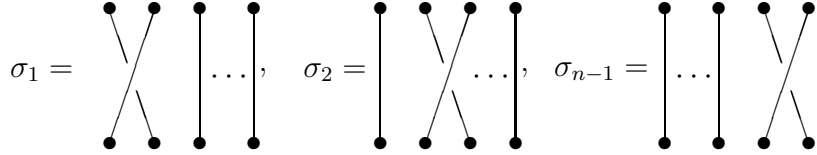
\begin{picture}(200,42)(-10,0)
\put(5,16){\makebox(0,0){\mbox{$\sigma_1=$}}}
\put(20,0){\makebox(0,0){\mbox{$\bullet$}}}
\put(30,0){\makebox(0,0){\mbox{$\bullet$}}}
\put(40,0){\makebox(0,0){\mbox{$\bullet$}}}
\put(52,0){\makebox(0,0){\mbox{$\bullet$}}}
\put(20,30){\makebox(0,0){\mbox{$\bullet$}}}
\put(30,30){\makebox(0,0){\mbox{$\bullet$}}}
\put(40,30){\makebox(0,0){\mbox{$\bullet$}}}
\put(52,30){\makebox(0,0){\mbox{$\bullet$}}}
\put(47,15){\makebox(0,0){\mbox{$\cdots$}}}
\put(20,0){\line(1,3){10}} \put(30,0){\line(-1,3){4}}
\put(20,30){\line(1,-3){4}} \put(40,0){\line(0,1){30}}
\put(52,0){\line(0,1){30}} \put(54,16){\makebox(0,0){\mbox{,}}}
\put(70,16){\makebox(0,0){\mbox{$\sigma_2=$}}}
\put(80,0){\makebox(0,0){\mbox{$\bullet$}}}
\put(90,0){\makebox(0,0){\mbox{$\bullet$}}}
\put(100,0){\makebox(0,0){\mbox{$\bullet$}}}
\put(110,0){\makebox(0,0){\mbox{$\bullet$}}}
\put(80,30){\makebox(0,0){\mbox{$\bullet$}}}
\put(90,30){\makebox(0,0){\mbox{$\bullet$}}}
\put(100,30){\makebox(0,0){\mbox{$\bullet$}}}
\put(110,30){\makebox(0,0){\mbox{$\bullet$}}}
\put(105,15){\makebox(0,0){\mbox{$\cdots$}}}
\put(80,0){\line(0,1){30}} \put(90,0){\line(1,3){10}}
\put(100,0){\line(-1,3){4}} \put(90,30){\line(1,-3){4}}
\put(110,0){\line(0,1){30}} \put(112,16){\makebox(0,0){\mbox{,}}}
\put(129,16){\makebox(0,0){\mbox{$\sigma_{n-1}=$}}}
\put(143,0){\makebox(0,0){\mbox{$\bullet$}}}
\put(155,0){\makebox(0,0){\mbox{$\bullet$}}}
\put(165,0){\makebox(0,0){\mbox{$\bullet$}}}
\put(175,0){\makebox(0,0){\mbox{$\bullet$}}}
\put(143,30){\makebox(0,0){\mbox{$\bullet$}}}
\put(150,15){\makebox(0,0){\mbox{$\cdots$}}}
\put(155,30){\makebox(0,0){\mbox{$\bullet$}}}
\put(165,30){\makebox(0,0){\mbox{$\bullet$}}}
\put(175,30){\makebox(0,0){\mbox{$\bullet$}}}
\put(143,0){\line(0,1){30}}
\put(165,0){\line(1,3){10}}
\put(175,0){\line(-1,3){4}}
\put(165,30){\line(1,-3){4}}
\put(155,0){\line(0,1){30}}
\end{picture}
}
\vskip 0.5cm
{\it Knot theory :} Alexander,
Jones,
HOMFLYPT,
Kauffman
polynomials.

Respectively: {\it Temperley-Lieb, Hecke, BMW} algebras.

{\it Geometry, physics} etc.

{\it Relation with the symmetric group} $S_n:\sigma_i^2=1$
$$\sigma_i^2=1\Rightarrow Sp\,(\rho(\sigma_i))\subseteq\{-1,1\}$$
$$
 Rep(S_n)\quad\quad
Rep(B_n) ?
$$
$$
\hat{S_n}=\{\text{\bf Young diagrams}\},\,\,
\text{\bf\,\,Plancherel \,\,measure\,\,on}\,\,\hat{S_n}.
$$
The {\it Young graph} explains how to {\it decompose the
restriction} $\rho\mid_{S_{n-1}}$ of the representation $\rho\in
\hat{S_n}$, etc.
$$\sigma_i\sigma_{i+1}\sigma_i=\sigma_{i+1}\sigma_i\sigma_{i+1}\Rightarrow
Sp\,(\rho(\sigma_i))=Sp\,(\rho(\sigma_{i+1})).$$
%
%
The {\it Hecke algebra} is defined by
$$H_n(q)=\langle\sigma_i)_{i=1}^{n-1}\mid ...\sigma_i^2=(q-1)\sigma_i+q
\rangle,\quad p_2(\sigma_i)=0,$$ hence
$Sp\,(\rho(\sigma_i))\subseteq\{-1,q\}$ and $H_n(q)\cong {\mathbb
C}[S_n]$.

\newpage
%

\newpage
$ 1.\,\, {\rm\bf Definition}\quad
B_3=\langle\sigma_1,\sigma_2\mid\sigma_1\sigma_2\sigma_1=\sigma_2\sigma_1\sigma_2
\rangle. $
\vskip 0.2cm
$2.\,\,{\bf Homomorphism}\quad\rho:B_3\mapsto {\rm SL}(2,{\mathbb
Z}),\,\,$
$$
 \sigma_1\mapsto\left(\begin{smallmatrix}
1&1\\
0&1
\end{smallmatrix}\right),\quad
\sigma_2\mapsto\left(\begin{smallmatrix}
1&0\\
-1&1
\end{smallmatrix}\right),\quad \sigma_2=(\sigma_1^{-1})^\sharp.
$$
\vskip 0.2cm
$ 3.\quad B_3/Z(B_3)\simeq{\rm PSL}(2,{\mathbb Z})\simeq {\mathbb
Z}_2\ast{\mathbb Z}_3. $
\vskip 0.2cm
4.\,\,{\bf P.~Humphries result, Pascal's triangle}
$$\sigma_1\mapsto
\sigma_1(1,n),\,\,\sigma_2\mapsto \sigma_2(1,n).$$
\vskip 0.2cm
5. \,\,{\bf Ferrand result} $ \Phi_n,\,\Psi_n\in{\rm
End}\,{\mathbb C}^n[X].$
\vskip 0.2cm
6.\,\, {\bf Tubo-Wenzl example}
$$\sigma_1,\mapsto\sigma_1(1,n)\Lambda_n,\quad
\sigma_2\mapsto \Lambda_n^\sharp\sigma_2(1,n),\quad
\Lambda_n\Lambda_n^\sharp=cI.$$

7.\,\, {\bf Tubo - Wenzl classifications of} $B_3-{\rm mod}$,
${\rm dim}V\leq 5.$
\vskip 0.2cm

8. \,\,{\bf Generalizations}
$$
\sigma_1\mapsto
\sigma_1^\Lambda(q,n):=\sigma_1(q,n)D_n(q)^\sharp\Lambda_n,\,\,
$$
$$
\sigma_2\mapsto\sigma_2^\Lambda(q,n):= \Lambda_n^\sharp
D_n(q)\sigma_2(q,n),
$$
$$
\text{where\,\,}\sigma_2(q,n)=(\sigma_1^{-1}(q^{-1},n))^\sharp,\,\,\Lambda_n={\rm
diag}(\lambda_r)_{r=0}^n,\,\,\Lambda_n\Lambda_n^\sharp=cI,
$$
$$
D_n(q)={\rm diag}(q_r)_{r=0}^n,\,\,
q_r=q^{\frac{(r-1)r}{2}},\,\,r,n\in{\mathbb N}.
$$
\vskip 0.2cm
9.\,\, {\bf The connection between $Rep(B_3)$ and}
$U_q(\mathfrak{sl}_2)$-mod.
\vskip 0.2cm
10.\,\, {\bf The Burau representation } $\rho_n:B_{n}\mapsto {\rm
GL}_{n}({\mathbb Z}[t,t^{-1}])$.
\vskip 0.2cm
11. {\bf Lowrence-Kramer representations}
\vskip 0.2cm
12.\,\, {\bf Generalization of 8 and 9 for} $B_n.$
\vskip 0.2cm
13.\,\, {\bf Formanek  classifications of} $B_n-{\rm mod}$, for
${\rm dim}V\leq n.$

\newpage
$ 1. \quad
B_3=\langle\sigma_1,\sigma_2\mid\sigma_1\sigma_2\sigma_1=\sigma_2\sigma_1\sigma_2
\rangle. $

$2.\quad\rho:B_3\mapsto {\rm SL}(2,{\mathbb Z}),\,\,$ $$
\sigma_1\mapsto\left(
\begin{array}{ccc}
1&1\\
0&1
\end{array}
\right),\quad \sigma_2\mapsto\left(
\begin{array}{ccc}
1&0\\
-1&1
\end{array}\right).
$$
$ 3.\quad B_3/Z(B_3)\simeq{\rm PSL}(2,{\mathbb Z})\simeq {\mathbb
Z}_2\ast{\mathbb Z}_3. $

 Hint: {\bf the Pascal triangle},
$\sigma_1\mapsto\sigma_2?\,\,$ $ \left(\begin{smallmatrix}
1&1\\
0&1
\end{smallmatrix}\right)^{-1}=\left(\begin{smallmatrix}
1&-1\\
0&1
\end{smallmatrix}\right).$
$$
\sigma_1(1,2):=\left(\begin{array}{ccc}
1&2&1\\
0&1&1\\
0&0&1
\end{array}\right),\quad
\sigma_1^{-1}(1,2)^\sharp:=\left(
\begin{array}{ccc}
1&0&0\\
-1&1&0\\
1&-2&1
\end{array}\right).
$$
Notations the {\bf central symmetry:}
$$
A^\sharp:=(A^t)^s,\quad A^\sharp=(a^\sharp_{ij}),\,\,
a^\sharp_{ij}=a_{n-i,n-j},
$$
$$
\sigma_1\mapsto \sigma_1(1,2),\quad \sigma_2\mapsto
\sigma_2(1,2):=\sigma_1^{-1}(1,2)^\sharp.
$$
\vskip 0.2cm
 4. {\bf P.~Humphries}, \cite{Hum00}
representations of $B_3$ in ${\mathbb C}^{n+1}$
\begin{equation}
\sigma_1\mapsto \sigma_1(1,n),\quad \sigma_2\mapsto
\sigma_2(1,n):=\sigma_1^{-1}(1,n)^\sharp.
\end{equation}
5. {\bf Ferrand  result}, \cite{Fer05}. $ \Phi_n,\,\Psi_n\in{\rm
End}\,
{\mathbb C}^n[X]:\Phi_n\Psi_n\Phi_n=\Psi_n\Phi_n\Psi_n.$
$$
(\Phi_np)(X):=p(X+1),\quad (\Psi_np)(X):=
(1-X)^np\left(X/(1-X)\right).
$$
6. {\bf Tubo-Wenzl example} \cite{TubWen01}: representations
$\sigma^\Lambda(1,n)$ of $B_3$ in ${\mathbb C}^{n+1}$
\begin{equation}
\sigma_1\mapsto \sigma_1(1,n)\Lambda_n,\quad \sigma_2\mapsto
\Lambda_n^\sharp\sigma_2(1,n),\,\,
\end{equation}
conditions on the complex diagonal matrix $\Lambda_n={\rm
diag}(\lambda_0,\lambda_1,...,\lambda_n)$ are the following:
\begin{equation}
\Lambda_n\Lambda_n^\sharp=cI,\,\,\,c\in {\mathbb C}.
\end{equation}

\newpage
{\bf 7.  Tubo - Wenzl classifications of $B_3-{\rm mod}$, ${\rm
dim}V\leq 5.$}

See \cite{TubWen01}. Let $V$ be a simple $B_3$ module of dimension
$n=2,3$. Then there exist a basis for $V$ for which $\sigma_1$ and
$\sigma_2$ act as follows ($\lambda=(\lambda_k)_k$) for  $n=2$ and
$n=3$
\begin{equation}
\label{TW2}
\sigma_1^{ \lambda}:=\left(\begin{smallmatrix}
\lambda_1&\lambda_1\\
0&\lambda_2&
\end{smallmatrix}\right)=
\left(\begin{smallmatrix}
\lambda_1&0\\
0&\lambda_2
\end{smallmatrix}\right)
\left(\begin{smallmatrix}
1&1\\
0&1
\end{smallmatrix}\right),\quad
\sigma_2^{ \lambda}:=\left(\begin{smallmatrix}
\lambda_2&0\\
-\lambda_2&\lambda_1
\end{smallmatrix}\right)=\left(\begin{smallmatrix}
1&0\\
-1&1
\end{smallmatrix}\right)\left(\begin{smallmatrix}
\lambda_2&0\\
0&\lambda_1
\end{smallmatrix}\right),\,\,
\end{equation}
\begin{equation}
\label{TW3}
\sigma_1\mapsto\sigma_1^{\lambda}=\left(\begin{smallmatrix}
\lambda_1&\lambda_1\lambda_3\lambda_2^{-1}+\lambda_2&\lambda_2\\
0&\lambda_2&\lambda_2\\
0&0&\lambda_3
\end{smallmatrix}\right),\quad
\sigma_2\mapsto\sigma_2^{ \lambda}:=\left(\begin{smallmatrix}
\lambda_3&0&0\\
-\lambda_2&\lambda_2&0\\
\lambda_2&-\lambda_1\lambda_3\lambda_2^{-1}-\lambda_2&\lambda_1
\end{smallmatrix}\right).
\end{equation}
Let us set $D=\sqrt{\lambda_2\lambda_3/\lambda_1\lambda_4}$. All
simple modules for $n=4$ are the following:
\begin{equation}
\label{TW41}
\sigma_1\mapsto\sigma_1^{\lambda}=\left(\begin{smallmatrix}
\lambda_1&(1+D^{-1}+D^{-2})\lambda_2&(1+D^{-1}+D^{-2})\lambda_3&\lambda_4\\
0&\lambda_2&(1+D^{-1})\lambda_3&\lambda_4\\
0&0&\lambda_3&\lambda_4\\
0&0&0&\lambda_4
\end{smallmatrix}\right),
\end{equation}
\begin{equation}
\label{TW42}
\sigma_2\mapsto\sigma_2^{\lambda}=\left(\begin{smallmatrix}
\lambda_4&0&0&0\\
-\lambda_3&\lambda_3&0&0\\
D\lambda_2&-(D+1)\lambda_2&\lambda_2&0\\
-D^3\lambda_1&(D^3+D^2+D)\lambda_1&-(D^2+D+1)\lambda_1&\lambda_1
\end{smallmatrix}\right).
\end{equation}
Let us set
$\gamma=(\lambda_1\lambda_2\lambda_3\lambda_4\lambda_5)^{1/5}$.
All simple modules for $n=5$ are the following:
\begin{equation}
\label{TW5} \sigma_1\mapsto\sigma_1^{
\lambda}=\left(\begin{smallmatrix} \lambda_1 &
(1+\frac{\gamma^2}{\lambda_2\lambda_4})(\lambda_2+\frac{\gamma^3}{\lambda_3\lambda_4})
&
(\frac{\gamma^2}{\lambda_3}+\lambda_3+\gamma)(1+\frac{\lambda_1\lambda_5}{\gamma^2})&
(1+\frac{\lambda_2\lambda_4}{\gamma^2})
(\lambda_3+\frac{\gamma^3}{\lambda_2\lambda_4})&\frac{\gamma^3}{\lambda_1\lambda_5}
\\
0&\lambda_2&\frac{\gamma^2}{\lambda_3}+\lambda_3+\gamma&
\frac{\gamma^3}{\lambda_1\lambda_5}+\lambda_3+\gamma&\frac{\gamma^3}{\lambda_1\lambda_5}
\\
0&0&\lambda_3&\frac{\gamma^3}{\lambda_1\lambda_5}+\lambda_3
&\frac{\gamma^3}{\lambda_1\lambda_5}\\
0&0&0&\lambda_4&\lambda_4\\
0&0&0&0&\lambda_5
\end{smallmatrix}\right).
\end{equation}
The formula for $\sigma_2^{ \lambda}$ was not given in
\cite{TubWen01}.
\newpage
{\bf 8. Equivalence of  Tuba-Wenzl's representations in the case
${\rm dim}\leq 5$ and our representations}.

General formulas for $1\leq n\leq 4$ gives us (we set
$q_r=q^{\frac{(r-1)r}{2}}$):
$$
\sigma_1\mapsto \sigma_1^\Lambda:=\sigma_1(q,n)\Lambda_n,\quad
\sigma_2\mapsto\sigma_2^\Lambda:= \Lambda_n^\sharp\sigma_2(q,n),
$$
$$
\Lambda_n\Lambda_n^\sharp=\lambda_0\lambda_n\Lambda_n(q),\quad
\Lambda_n(q)=q_n^{-1}D_n(q)D_n^\sharp(q),\,\,D_n(q)={\rm
diag}(q_r)_{r=0}^n,
$$
\begin{equation}
\label{cond}
\lambda_r\lambda_{n-r}=\lambda_0\lambda_nq^{-(n-r)r},\quad 0\leq
r\leq n.
\end{equation}
Let $n=1$ we have
$$
\sigma_1^\Lambda
=\left(\begin{smallmatrix}
1&1\\
0&1
\end{smallmatrix}\right)\Lambda_1,\quad
\sigma_2^\Lambda
=\Lambda_1^\sharp\left(\begin{smallmatrix}
1&0\\
-1&1
\end{smallmatrix}\right),\quad
 \Lambda_1=\left(\begin{smallmatrix}
\lambda_0&0\\
0&\lambda_1
\end{smallmatrix}\right).
$$
Let $n=2$, conditions (\ref{cond}) gives us $ \Lambda_2={\rm
diag}(\lambda_r)_{r=0}^3$
$${\rm
diag}(\lambda_0\lambda_2,\lambda_1^2,\lambda_0\lambda_2)=\lambda_0\lambda_2{\rm
diag}(1,q^{-1},1),\text{\,\,\, so\,\,\,}
q^{-1}=\lambda_1^2/\lambda_0\lambda_2.
$$
$$
\sigma_1^\Lambda(q,2)=\left(\begin{smallmatrix}
1&1+q&1\\
0&1&1\\
0&0&1
\end{smallmatrix}\right)\Lambda_2,\quad
\sigma_2^\Lambda (q,2) =
\Lambda_2^\sharp\left(\begin{smallmatrix}
1&0&0\\
-1&1&0\\
 q^{-1}&-(1+q^{-1})&1
\end{smallmatrix}\right).
$$
For $n=3$ conditions (\ref{cond}) gives us
$q^{-2}=\lambda_1\lambda_2/\lambda_0\lambda_3$ for $r=1$.
$$
\sigma_1(q,3)=\left(\begin{smallmatrix}
1&1+q+q^2&1+q+q^2&1\\
0&1&1+q&1\\
0&0&1&1\\
0&0&0&1
\end{smallmatrix}\right),\quad
\Lambda=\left(\begin{smallmatrix}
\lambda_0&0&0&0\\
0&\lambda_1&0&0\\
0&0&\lambda_2&0\\
0&0&0&\lambda_3
\end{smallmatrix}\right),
$$
$$
\sigma_2(q,3) = \left(\begin{smallmatrix}
 1    &0          &0              &0\\
-1    &1         &0                &0\\
q^{-1}&-(1+q^{-1})&1                 &0\\
-q^{-3}&q^{-1}(1+q^{-1}+q^{-2})&-(1+q^{-1}+q^{-2})&1
\end{smallmatrix}\right).
$$
For $n=4$ conditions (\ref{cond}) gives us
$q^{-3}=\lambda_1\lambda_3/\lambda_0\lambda_4$ for $r=1$ and
$q^{-4}=\lambda_2^2/\lambda_0\lambda_4$ for $r=2$.
$$
\sigma_1(q)=\left(\begin{smallmatrix}
1&(1+q)(1+q^2)&(1+q^2)(1+q+q^2)&(1+q)(1+q^2)&1\\
0&1&1+q+q^2&1+q+q^2&1\\
0&0&1&1+q&1\\
0&0&0&1&1\\
0&0&0&0&1
\end{smallmatrix}\right),\,\,
\Lambda=\left(\begin{smallmatrix}
\lambda_0&0&0&0&0\\
0&\lambda_1&0&0&0\\
0&0&\lambda_2&0&0\\
0&0&0&\lambda_3&0\\
0&0&0&0&\lambda_4
\end{smallmatrix}\right),
$$
$\sigma_2(q,4)=(\sigma_1^{-1}(q^{-1},4))^\sharp$.

\newpage
$$
\sigma_1\mapsto \sigma_1(1,n)\Lambda_n,\quad \sigma_2\mapsto
\Lambda_n^\sharp\sigma_2(1,n), (2)
$$
$$\Lambda_n={\rm diag}(\lambda_r)_{r=0}^n,\quad
\Lambda\Lambda^\sharp=cI,\,\,\,c\in {\mathbb C,}\quad (3)
$$
{\bf 8. Generalization} of (2) for $q\not=1,$  with the condition
(3)
\begin{equation}
\label{rep} \sigma_1\mapsto \sigma_1^\Lambda(q,n):=
\sigma_1(q,n)D_n^\sharp(q)\Lambda_n,\,\, \sigma_2\mapsto
\sigma_2^\Lambda(q,n):= \Lambda_n^\sharp D_n(q)\sigma_2(q,n),
\end{equation}
\begin{equation}
\sigma_2(q,n):=\sigma_1^{-1}(q^{-1},n)^\sharp,\,\, D_n(q)={\rm
diag}(q_r)_{r=0}^n,\,\,q_r=q^{\frac{(r-1)r}{2}},
\end{equation}
where $q-$ {\it binomial coefficients} or  {\it Gaussian
polynomials} are defined as follows
\begin{equation}
\label{GP} \left(\begin{smallmatrix} n\\
k
\end{smallmatrix}\right)_q:=\frac{(n)!_q}{(k)!_q(n-k)!_q},\quad
\left[\begin{smallmatrix} n\\
k
\end{smallmatrix}\right]_q:=\frac{[n]!_q}{[k]!_q[n-k]!_q}
\end{equation}
corresponding to two forms of $q-${\it natural numbers}, defined
by
\begin{equation}
\label{(),[]_q}
(n)_q:=\frac{q^n-1}{q-1},\quad
[n]_q:=\frac{q^n-q^{-1}}{q-q^{-1}}.
\end{equation}

\begin{thm} {\rm \cite{KosAlb07q}} The formulas (\ref{rep}) $\sigma_1\mapsto
\sigma_1^\Lambda(q,n),\,\,\sigma_2\mapsto\sigma_2^\Lambda(q,n)$
give the representation of $B_3.$
\end{thm}
\begin{thm} {\rm \cite{KosAlb07q}} The representation $\sigma^\Lambda(q,n)$ defined
by (\ref{rep}) generalize the Tubo-Wenzl representations for
arbitrary $n\in {\mathbb N}.$
\end{thm}

{\bf Definition}. {\it We say that the representation is {\bf
subspace irreducible} or {\bf ireducible} (resp. {\bf operator
irreducible}) when there no nontrivial invariant close {\bf
subspaces} for all operators of the representation (resp.  there
no nontrivial bounded {\bf operators} commuting with all operators
of the representation).}

Let us define for $n,r,q,\lambda$ such that $n\in{\mathbb
N},\,\,0\leq r\leq n,\,\,\lambda\in{\mathbb
C}^{n+1},\,\,q\in{\mathbb C}$ the following operators
%
%
\begin{equation}
\label{F(nu)}
 F_{r,n}(q,\lambda)=\exp_{(q)}\left(\sum_{k=0}^{n-1}
(k+1)_qE_{kk+1}\right)-
q_{n-r}\lambda_r(D_n(q)\Lambda_n^\sharp)^{-1},
\end{equation}
where $\exp_{(q)} X=\sum_{m=0}^\infty X^m/(m)!_q$.
For the matrix $C\in{\rm Mat}(n+1,{\mathbb C})$ we denote by
$$
M^{i_1i_2...i_r}_{j_1j_2...j_r}(C),\,\,({\rm resp.\,\,}
A^{i_1i_2...i_r}_{j_1j_2...j_r}(C)),\,\, 0\leq i_1<...<i_r\leq
n,\,\, 0\leq j_1<...<j_r\leq n
$$
its  minors (resp. the cofactors) with $i_1,i_2,...,i_r$ rows and
$j_1,j_2,...,j_r$ columns.
\begin{thm}
\label{t.IrrB_3} {\rm \cite{KosAlb07q}} The representation of the
group $B_3$ defined by (10)
have the following properties:\\
1) for $q=1,\,\,\Lambda_n=1$, it is  {\rm subspace irreducible} in arbitrary dimension $n\in{\mathbb N }$;\\
2) for $q\not=1,\,\,\Lambda_n={\rm diag}(\lambda_k)_{k=0}^n\not=1$
it is {\rm operator irreducible} if and only if for any $0\leq
r\leq \left[\frac{n}{2}\right]$ there exists $0\leq
i_0<i_i<...<i_r\leq n $ such that
\begin{equation}
\label{M^i_j...=0} M^{i_0i_i...i_{n-r-1}}_{r+1r+2...n} (
F_{r,n}^s(q,\lambda))\not=0;
\end{equation}
3) for $q\not=1,\,\,\Lambda_n=1$ it is  {\it subspace irreducible}
if and only if $(n)_q\not=0$.\\
 The representation has $[\frac{n+1}{2}]+1$ free parameters.
\end{thm}
{\bf 9. The connection between $Rep(B_3)$ and
$U_q(\mathfrak{sl}_2)$-mod.}

The algebra $U(\mathfrak{sl}_2)$ is the associative algebra
generated by three generators $X,\,Y,\,H$ with the relations (7).
\begin{equation}\label{sl_2}
[H,X]=2X,\,\,[H,Y]=-2Y,\,\,[X,Y]=H,
\end{equation}
$$
X= \left(\begin{smallmatrix}
0&1\\
0&0
\end{smallmatrix}\right),\quad Y=
\left(\begin{smallmatrix}
0&0\\
1&0
\end{smallmatrix}\right),\quad H=\left(\begin{smallmatrix}
1&0\\
0&-1
\end{smallmatrix}\right) \text{\quad in\quad} \mathfrak{sl}_2.
$$
$U_q(\mathfrak{sl}_2)$ is the algebra generated by four variables
$E,\,F,\,K,\,K^{-1}$ with the relations
\begin{equation}
\label{U_q(1)} KK^{-1}=K^{-1}K=1,
\end{equation}
\begin{equation}
\label{U_q(2)} KEK^{-1}=q^2E,\quad KFK^{-1}=q^{-2}F,
\end{equation}
\begin{equation}
\label{U_q(3)}
[E,F]=\frac{K-K^{-1}}{q-q^{-1}}=\frac{q^H-q^{-H}}{q-q^{-1}}.
\end{equation}
Comultiplication $\Delta$, counit $\varepsilon$ and antipod $S$
are as follows:
$$
\Delta(E)=E\otimes K+1\otimes E, \quad \Delta(F)=F\otimes 1+
K^{-1}\otimes F,\quad \Delta(K)=K\otimes K,
$$
$$
S(K)=K^{-1},\,\,S(E)=-EK^{-1},\,\,S(F)=-KF,
$$
$$
\,\,\varepsilon(K)=1,\,\, \varepsilon(E)=\varepsilon(F)=0.
$$
 All finite-dimensional $U-$module $V$ being the highest weight module of highest weight
$\lambda$ are of the following form (see Kassel,
\cite[TheoremV.4.4.]{Kas95})
$$ \rho(n)(X)=\left(\begin{smallmatrix}
0&n&0&...&0\\
0&0&n-1&...&0\\
&&&...&\\
0&0&0&...&1\\
0&0&0&...&0
\end{smallmatrix}\right),\quad
 \rho(n)(Y)=\left(\begin{smallmatrix}
0&0  &...&0  &0\\
1&0  &...&0  &0\\
0&2&...&0  &0\\
 &   &...&   &\\
0&0  &...&n&0
\end{smallmatrix}\right),\,
$$
$$
\rho(n)(H)=\left(\begin{smallmatrix}
n&0      &...&0       &0\\
0  &n-2&...&0       &0\\
   &       &...&        &\\
   &       &...&-n+2&0\\
0  &0      &...&0       &-n\\
\end{smallmatrix}\right).
$$
where $\lambda ={\rm dim}(V)-1\in{\mathbb N}.$

 All finite-dimensional $U_q-$module $V$ being the highest weight module of highest weight
$\lambda$ are of the following form (see Kassel,
\cite[Theorem VI.3.5.]{Kas95})
$$ \rho_{\varepsilon,n}(E)=\varepsilon\left(\begin{smallmatrix}
0&[n]&0&...&0\\
0&0&[n-1]&...&0\\
&&&...&\\
0&0&0&...&1\\
0&0&0&...&0
\end{smallmatrix}\right),\quad
 \rho_{\varepsilon,n}(F)=\left(\begin{smallmatrix}
0&0  &...&0  &0\\
1&0  &...&0  &0\\
0&[2]&...&0  &0\\
 &   &...&   &\\
0&0  &...&[n]&0
\end{smallmatrix}\right),\,
$$
$$
\rho_{\varepsilon,n}(K)=\varepsilon\left(\begin{smallmatrix}
q^n&0      &...&0       &0\\
0  &q^{n-2}&...&0       &0\\
   &       &...&        &\\
   &       &...&q^{-n+2}&0\\
0  &0      &...&0       &q^{-n}\\
\end{smallmatrix}\right),
$$
where $\varepsilon=\pm 1,\,\,\lambda=\varepsilon q^n$ and
$n\in{\mathbb N}$.

 {\bf The main observation is the following:}
$$
\left(\begin{smallmatrix}
1&2&1\\
0&1&1\\
0&0&1
\end{smallmatrix}\right)=\exp\left(\begin{smallmatrix}
0&2&0\\
0&0&1\\
0&0&0
\end{smallmatrix}\right),\quad \left(\begin{smallmatrix}
1&(2)_q&1\\
0&1&(1)_q\\
0&0&1
\end{smallmatrix}\right)=\exp_{(q)}\left(\begin{smallmatrix}
0&(2)_q&0\\
0&0&(1)_q\\
0&0&0
\end{smallmatrix}\right),
$$
where
$$
\left(\begin{smallmatrix}
0&(2)_{q^2}&0\\
0&0&(1)_{q^2}\\
0&0&0
\end{smallmatrix}\right)=\left(\begin{smallmatrix}
0&[2]_q&0\\
0&0&[1]_q\\
0&0&0
\end{smallmatrix}\right)\left(\begin{smallmatrix}
q^2&0&0\\
0&q&0\\
0&0&1
\end{smallmatrix}\right),\,\,
\exp_{(q)} X:=\sum_{m=0}^\infty\frac{1}{(m)!_q}X^m.
$$
\begin{thm} For $q=1$ holds
\begin{equation}
\label{exp_1r(X)} \sigma_1(1,n)=\exp\left( \rho(n)(X)\right),\quad
\sigma_2(1,n)=\exp \left(\rho(n)(-Y)\right).
\end{equation}
\end{thm}
\begin{thm}
For $q\not=1$ we have
\begin{equation}
\label{exp_qr(E)}
\sigma_1(q^2,n)D_n^\sharp(q^2)=\exp_{(q^2)}\left(q^{n/2}\rho_{1,n}(EK^{1/2}
)\right) D_n^\sharp(q^2),
\end{equation}
\begin{equation}
\label{exp_qr(F)} \quad
D_n(q^2)\sigma_2(q^2,n)=\exp_{(q^2)}\left(-q^{n/2}\rho_{1,n}(FK^{-1/2}
)\right)D_n(q^2).
\end{equation}
\end{thm}

{\bf Proof.} The two
forms of $q-$natural numbers are connected as follows (see Kassel,
\cite{Kas95})
\begin{equation}
\label{[n]=(n))}
 [n]= q^{-(n-1)}(n)_{q^2},\quad
[n]!=q^{-(n-1)n/2}(n)!_{q^2}
\end{equation}
$$
\left(\begin{smallmatrix}
0&(n)&0&...&0\\
0&0&(n-1)&...&0\\
&&&...&\\
0&0&0&...&(1)\\
0&0&0&...&0
\end{smallmatrix}\right)=\left(\begin{smallmatrix}
0&[n]&0&...&0\\
0&0&[n-1]&...&0\\
&&&...&\\
0&0&0&...&[1]\\
0&0&0&...&0
\end{smallmatrix}\right){\rm diag}(q^n,q^{n-1},...,1)
$$
$=q^{n/2}\rho_{1,n}(EK^{1/2}
), $ and
$$
\left(\begin{smallmatrix}
0&0  &...&0  &0\\
(1)&0  &...&0  &0\\
0&(2)&...&0  &0\\
 &   &...&   &\\
0&0  &...&(n)&0
\end{smallmatrix}\right)=
\left(\begin{smallmatrix}
0&0  &...&0  &0\\
[1]&0  &...&0  &0\\
0&[2]&...&0  &0\\
 &   &...&   &\\
0&0  &...&[n]&0
\end{smallmatrix}\right){\rm diag}(1,q,...,q^{n-1},q^{n})
$$
$= q^{n/2}\rho_{1,n}(FK^{-1/2}
),$
since
$$
{\rm diag}(1,q,...,q^{n-1},q^{n})
=q^{n/2}\rho_{1,n}(K^{-1/2}
)
$$
and
$$
{\rm diag}(q^n,q^{n-1},...,1)=
q^{n/2}\rho_{1,n}(K^{1/2}
) .
$$
Al last we conclude that $$ \left(\begin{smallmatrix}
0&(n)&0&...&0\\
0&0&(n-1)&...&0\\
&&&...&\\
0&0&0&...&(1)\\
0&0&0&...&0
\end{smallmatrix}\right)=q^{n/2}\rho_{1,n}(EK^{1/2}
),\,\,
$$
$$
\left(\begin{smallmatrix}
0&0  &...&0  &0\\
(1)&0  &...&0  &0\\
0&(2)&...&0  &0\\
 &   &...&   &\\
0&0  &...&(n)&0
\end{smallmatrix}\right)=q^{n/2}\rho_{1,n}(FK^{-1/2}
).
$$
Further we observe that
$$
X\otimes I+I\otimes X\mid_{S^{2}({\mathbb C}^2)}=
\left(\begin{smallmatrix}
0&1\\
0&0
\end{smallmatrix}\right)\otimes I+ I\otimes\left(\begin{smallmatrix}
0&1\\
0&0
\end{smallmatrix}\right)\mid_{S^{2}({\mathbb
C}^2)}=\left(\begin{smallmatrix}
0&2&0\\
0&0&1\\
0&0&0
\end{smallmatrix}\right)
$$
$$
\Delta\rho(1)(X)\mid_{S^{2}({\mathbb C}^2)}=\rho(2)(X),
$$
$$
(I+X)\otimes(I+X)=\exp(\Delta (X)),\quad \sigma_1(1,1)\otimes
\sigma_1(1,1)\mid_{S^{2}({\mathbb C}^2)}=\sigma (1,2).
$$
\begin{lem}
\label{l.Sym^n(1)=(n)}
 We have for $q\not=1$
\begin{equation}
\label{Sym^n(1)=(n)} \rho_{1,n}=
\Delta^{n-1}\rho_{1,1}\mid_{S^{n,q}({\mathbb C}^2)},
\end{equation}
where $ S^{n,q}({\mathbb C}^2)$ is $q-$symmetric tensor power of
${\mathbb C}^2$.
\end{lem}
{\bf Proof.} For $n=1$ we have the following operators
$$
\rho_{1,1}(E)=\left(\begin{smallmatrix}
0&1\\
0&0
\end{smallmatrix}\right)
,\quad
 \rho_{1,1}(F)=\left(\begin{smallmatrix}
0&0\\
1&0
\end{smallmatrix}\right)
,\quad\rho_{1,1}(K)= \left(\begin{smallmatrix}
q&0\\
0&q^{-1}
\end{smallmatrix}\right)= q^{H}.
$$
For $n=2$ we get
$$
\rho_{1,2}(E)=
\left(\begin{smallmatrix}
0&[2]&0\\
0&0&[1]\\
0&0&0\\
\end{smallmatrix}\right),\,\,
\rho_{1,
2}(F) =\left(\begin{smallmatrix}
0&0&0\\
[1]&0&0\\
0&[2]&0\\
\end{smallmatrix}\right),\,\,
\rho_{1,
2}(K) =
\left(\begin{smallmatrix}
q^2&0&0\\
0&1&0\\
0&0&q^{-2}\\
\end{smallmatrix}\right)
$$
We have $\Delta (\rho_{1,1}(E))=$
$$ \rho_{1,1}(E)\otimes
\rho_{1,1}(K)+1\otimes \rho_{1,1}(E)=
\left(\begin{smallmatrix}
0&1\\
0&0
\end{smallmatrix}\right)\otimes \left(\begin{smallmatrix}
q&0\\
0&q^{-1}
\end{smallmatrix}\right)+\left(\begin{smallmatrix}
1&0\\
0&1
\end{smallmatrix}\right)\otimes \left(\begin{smallmatrix}
0&1\\
0&0
\end{smallmatrix}\right)
$$
$$
=\left(\begin{smallmatrix}
0&0&q&0\\
0&0&0&q^{-1}\\
0&0&0&0\\
0&0&0&0
\end{smallmatrix}\right)+
\left(\begin{smallmatrix}
0&1&0&0\\
0&0&0&0\\
0&0&0&1\\
0&0&0&0
\end{smallmatrix}\right)=
 \left(\begin{smallmatrix}
0&1&q&0\\
0&0&0&q^{-1}\\
0&0&0&1\\
0&0&0&0
\end{smallmatrix}\right).
$$
Further $\Delta (\rho_{1,1}(F))=$
$$
 \rho_{1,1}(F)\otimes 1+\rho_{1,1}(K^{-1})\otimes
\rho_{1,1}(F)=
 \left(\begin{smallmatrix}
0&0\\
1&0
\end{smallmatrix}\right)\otimes
\left(\begin{smallmatrix}
1&0\\
0&1
\end{smallmatrix}\right)+
\left(\begin{smallmatrix}
q^{-1}&0\\
0&q
\end{smallmatrix}\right)
\otimes \left(\begin{smallmatrix}
0&0\\
1&0
\end{smallmatrix}\right)
$$
$$
=\left(\begin{smallmatrix}
0&1&0&0\\
0&0&0&0\\
1&0&0&1\\
0&1&0&0
\end{smallmatrix}\right)+
\left(\begin{smallmatrix}
0&0&0&0\\
q^{-1}&0&0&0\\
0&0&0&0\\
0&0&q&0
\end{smallmatrix}\right)=
\left(\begin{smallmatrix}
0     &0&0&0\\
q^{-1}&0&0&0\\
1     &0&0&0\\
0     &1&q&0
\end{smallmatrix}\right)
$$
and
$$
\Delta (\rho_{1,1}(K))=\rho_{1,1}(K)\otimes \rho_{1,1}(K)=
\left(\begin{smallmatrix}
q&0\\
0&q^{-1}
\end{smallmatrix}\right)\otimes\left(\begin{smallmatrix}
q&0\\
0&q^{-1}
\end{smallmatrix}\right)=\left(\begin{smallmatrix}
q^2&0&0&0\\
0  &1&0&0\\
0  &0&1&0\\
0  &0&0&q^{-2}
\end{smallmatrix}\right).
$$
In the {\it $q-$symmetric basis} of the submodule
$S^{2,q}({\mathbb C}^2)$ of the module ${\mathbb C}^2\otimes
{\mathbb C}^2$
$$
e_{00}^{s,q}=e_0\otimes e_0,\quad e_{01}^{s,q}=q^{-1}e_0\otimes
e_1+e_1\otimes e_0,\quad e_{11}^{s,q}=e_1\otimes e_1
$$
the operator $\Delta (\rho_{1,1}(E))$ has the following form:
$$
 \Delta
(\rho_{1,1}(E))\mid_{S^{2,q}({\mathbb
C}^2)}=\left(\begin{smallmatrix}
0&[2]&0\\
0&0&[1]\\
0&0&0\\
\end{smallmatrix}\right).
$$
The basis in the space ${\mathbb C}^2\otimes {\mathbb C}^2$ is
generated by vectors $e_{kn},\,0\leq k,n\leq 1$ where
$e_{kn}=e_k\otimes e_n$. Operator $\Delta (\rho_{1,1}(E))$ acts as
follows $ e_{00}\mapsto 0,\,\,e_{01}\mapsto e_{00},\,\,
e_{10}\mapsto qe_{00} ,\,\,e_{11}\mapsto\,\,q^{-1}e_{01}+e_{10}$,
hence $e_{00}^{s,q}\mapsto 0,\,\,$
$$
e_{01}^{s,q}=q^{-1}e_{01}+e_{10}\mapsto(q+q^{-1})e_{00}=[2]e_{00}^{s,q}
,\,\, e_{11}^{s,q}\mapsto\,\,q^{-1}e_{01}+e_{10}=e_{01}^{s,q}.
$$
Similarly we get
$$
\Delta (\rho_{1,1}(F))\mid_{S^{2,q}({\mathbb
C}^2)}=\left(\begin{smallmatrix}
0&0&0\\
[1]&0&0\\
0&[2]&0\\
\end{smallmatrix}\right),\quad
\Delta (\rho_{1,1}(K))\mid_{S^{2,q}({\mathbb
C}^2)}=\left(\begin{smallmatrix}
q^2&0&0\\
0&1&0\\
0&0&q^{-2}\\
\end{smallmatrix}\right).
$$
hence (\ref{Sym^n(1)=(n)}) holds for $n=2$. For $n>2$ the proof is
similar.
%
\vskip 0.5cm
 {\bf 10. The Burau representation} $\rho:B_{n}\mapsto
{\rm GL}_n({\mathbb Z}[t,t^{-1}])$  is defined
for  a non-zero complex number $t$ by
$$
\sigma_i\mapsto\beta_i=I_{i-1}\oplus\left(\begin{smallmatrix}
1-t&t\\
1&0\\
\end{smallmatrix}\right)\oplus I_{n-i-1}
$$
 where $1-t$ is the $(i,i)$ entry.
Representation $\rho$ splits into 1-dimesional and
$n\!-\!1-$dimensional irreducible representations, known as {\it
reduced Burau representation} $\overline{\rho}:B_{n}\mapsto {\rm
GL}_{n-1}({\mathbb Z}[t,t^{-1}])$
$$
\sigma_1\mapsto b_1=\left(\begin{smallmatrix}
-t&0\\
-1&1\\
\end{smallmatrix}\right)\oplus I_{n-3},\quad
\sigma_{n-1}\mapsto b_{n-1}=I_{n-3}\oplus
\left(\begin{smallmatrix}
1&-t\\
0&-t\\
\end{smallmatrix}\right),
$$
$$
\sigma_i\mapsto b_i=I_{i-2}\oplus\left(\begin{smallmatrix}
1&-t&0\\
0&-t&0\\
0&-1&1\\
\end{smallmatrix}\right)\oplus I_{n-i-2},\,\,2\leq i\leq n-2.
$$
\vskip 0.5cm
{\bf Problem.} Whether the reduced Burau representation
$\overline{\rho}:B_n\mapsto {\rm GL}_{n-1}({\mathbb Z}[t,t^{-1}])$
is {\it faithful} ?

YES for $n=3$ (Birman
 \cite{Bir74}).
NO for $n\geq 9$ Moody
\cite{Mod91} Long and Paton \cite{LonPat93},
Bigelow \cite{Big99}
improved further for $n\geq 5$.
%
\vskip 0.5cm
{\bf Open problem:} Whether the reduced Burau representation of
$B_4\mapsto {\rm GL}_3({\mathbb Z}[t,t^{-1}])$
$$
b_1=\left(\begin{smallmatrix}
-t&0&0\\
-1&1&0\\
0&0&1\\
\end{smallmatrix}\right),\,\,
b_2=\left(\begin{smallmatrix}
1&-t&0\\
0&-t&0\\
0&-1&1\\
\end{smallmatrix}\right),\,\,
b_3=\left(\begin{smallmatrix}
1&0&0\\
0&1&-t\\
0&0&-t\\
\end{smallmatrix}\right)
$$
is {\bf faithful}
\vskip 0.5cm
{\bf 11. Lowrence-Kramer representations,
\cite{Kra00}}
$$\lambda : B_n\mapsto
{\rm GL}_m({\mathbb Z}[t^{\pm 1},q^{\pm 1}]),\quad m=n(n-1)/2.$$
The basis in the space ${\mathbb C}^{n(n-1)/2}$ is $x_{ik},\,1\leq
i<k\leq n.$

{\bf Faithfulness for all} $n$, Bigelow
\cite{Big01}, Kramer
\cite{Kra02} $\Rightarrow B_n$ {\bf is a linear group for all}
$n.$
$$
\left.\begin{array}{rll}
\sigma_kx_{k,k+1}=&tq^2x_{k,k+1}&\\
\sigma_kx_{ik}=&(1-q)x_{ik}+qx_{i,k+1}&\text{\,for\,\,}i<k\\
\sigma_kx_{i,k+1}=&x_{ik}+tq^{k-i+1}(q-1)x_{k,k+1}&\text{\,for\,\,}i<k\\
\sigma_kx_{kj}=&tq(q-1)x_{k,k+1}+qx_{k+1,j}&\text{\,for\,\,}k+1<j\\
\sigma_kx_{k+1,j}=&x_{kj}+(1-q)x_{k+1,j}&\text{\,for\,\,}k+1<j\\
\sigma_kx_{ij}=&x_{ij}&\text{\,for\,\,}i<j<k\text{\,or\,}k+1<i<j\\
\sigma_kx_{ij}=&x_{ij}+tq^{k-i}(q-1)^2x_{k,k+1}&\text{\,for\,\,}i<k<k+1<j\\
\end{array}\right.
$$
\vskip 0.5cm
{\bf 12. Generalization of 8 and 9 for $B_n$.} For $n=4$ and
$t=-1$ we have $\overline{\rho}_4:B_4\mapsto {\rm SL}(3,{\mathbb
Z})$
$$
b_1=\left(\begin{smallmatrix}
1&0&0\\
-1&1&0\\
0&0&1\\
\end{smallmatrix}\right),\,\,
b_2=\left(\begin{smallmatrix}
1&1&0\\
0&1&0\\
0&-1&1\\
\end{smallmatrix}\right),\,\,
b_3=\left(\begin{smallmatrix}
1&0&0\\
0&1&1\\
0&0&1\\
\end{smallmatrix}\right).
$$
$$
b_1=\exp(-F_1),\,b_2=\exp(E_1-F_2),b_3=\exp(E_2).
$$
We can show that the symmetric powers $b_i\otimes b_i\mid_S$ are
the following
$$
b_1\otimes b_1\mid_S=\left(\begin{smallmatrix}
1&0&0&0&0&0\\
-1&1&0&0&0&0\\
1&-2&1&0&0&0\\
0&0&0&1&0&0\\
0&0&0&-1&1&0\\
0&0&0&0&0&1\\
\end{smallmatrix}\right),\,
b_2\otimes b_2\mid_S=\left(\begin{smallmatrix}
1&2&1&0&0&0\\
0&1&1&0&0&0\\
0&0&1&0&0&0\\
0&-1&0&1&1&0\\
0&0&-2&0&1&0\\
0&0&1&0&-1&1\\
\end{smallmatrix}\right),
$$
$$
b_3\otimes b_3\mid_S= \left(\begin{smallmatrix}
1&0&0&0&0&0\\
0&1&0&1&0&0\\
0&0&1&0&2&1\\
0&0&0&1&0&0\\
0&0&0&0&1&1\\
0&0&0&0&0&1\\
\end{smallmatrix}\right).
$$
We have for $n=5$ and $t=-1$ $b^{(5)}:B_5\mapsto {\rm
SL}(4,{\mathbb Z})$
$$
b_1=\left(\begin{smallmatrix}
1&0&0&0\\
-1&1&0&0\\
0&0&1&0\\
0&0&0&1\\
\end{smallmatrix}\right),\,\,
b_2=\left(\begin{smallmatrix}
1&1&0&0\\
0&1&0&0\\
0&-1&1&0\\
0&0&0&1\\
\end{smallmatrix}\right),\,\,
b_3=\left(\begin{smallmatrix}
1&0&0&0\\
0&1&1&0\\
0&0&1&0\\
0&0&-1&1\\
\end{smallmatrix}\right),\,\,
b_4=\left(\begin{smallmatrix}
1&0&0&0\\
0&1&0&0&\\
0&0&1&1\\
0&0&0&1\\
\end{smallmatrix}\right).
$$
Let  $\overline{\rho}:B_{n}\mapsto {\rm SL}_{n-1}({\mathbb Z})$ be
the {\it reduced Burrau representation} for $t=-1$.

The {\it quantum group} $U_q(\mathfrak{sl}_{n-1})$ is the algebra
generated by $4(n-1)$ variables $E_i,\,F_i,\,K_i,\,K^{-1}_i$ with
relations as (17)--(19). Let
$$
\rho_m:U_q(\mathfrak{sl}_{n-1})\mapsto {\rm End}({\mathbb C}^m)
$$
be the highest weight $U_q(\mathfrak{sl}_{n-1})-$module. {\it Then
$$
\sigma_1\mapsto\exp(-\rho_m(F_1)),\,\sigma_k\mapsto\exp(\rho_m(E_{k-1}-F_k)),\,\,
\sigma_n\mapsto\exp(\rho_m(E_{n-1})).
$$
gives the representation of $B_{n}$ for} $q=1$ (see (20)).

For $q\not=1$ we can obtain formulas similar to (21)--(22).
\vskip 0.5cm
{\bf 13.  Formanek  classifications of $B_n-{\rm mod}$, for ${\rm
dim}V\leq n.$}
\vskip 0.2cm
 In \cite{ForWSV03} E.Formanek et al.
gave the {\it complete classification} of all {\it simple
representations of $B_n$} for {\it dimension} $\leq n$.
%
\vskip 0.5cm
{\bf Acknowledgements.} {The author would like to thank the
Max-Planck-Institute of Mathematics and
the Institute of Applied Mathematics, University of Bonn
for the hospitality. The partial financial support by the DFG
project 436 UKR 113/87 is gratefully acknowledged.}


\end{document}